\documentclass{amsart}
\UseRawInputEncoding 
\usepackage{amsfonts}
\usepackage{amsmath}
\usepackage{amssymb}
\usepackage{graphicx}%
\providecommand{\U}[1]{\protect\rule{.1in}{.1in}}

\theoremstyle{plain}

\numberwithin{equation}{section}

\begin{document}

\title[ ]{Extension properties for monotone and sublinear operators}
\author{Sorin G. Gal}
\address{Department of Mathematics and Computer Science\\
University of Oradea\\
University\ Street No. 1, Oradea, 410087, Romania}
\email{galso@uoradea.ro, galsorin23@gmail.com}

\date{}


\begin{abstract}
In this paper we obtain several extension properties for  monotone and sublinear operators. The results obtained generalize those known for positive and linear operators.
\end{abstract}
\maketitle

\section{Introduction}

Recently, the classical Korovkin's approximation theory in function spaces for positive linear operators was extented to the
framework of sublinear and monotone operators by the papers Gal-Niculescu \cite{Gal-Nic-Med}--\cite{submit}.

This naturally suggested us that other results too for positive linear operators could be extended to the more general class of monotone sublinear
operators.

The aim of the present paper is to extend the study of monotone sublinear operations
concerning their extension properties. Although the proofs mostly follow the standard proofs in the linear case, but adapted to the sublinear one,
in my opinion the results obtained are new and good to be known.

In Section 2 we present some preliminaries on monotone sublinear operators.
Section 3 deals with several extension results on monotone sublinear operators, including those of Kantorovich-type and of Hahn-Banach-type.

\section{Preliminaries on sublinear monotone operators}

Let $E$ and $F$ be two ordered vector spaces. The following concepts are well-known.

{\bf Definition 2.1.} (i) The operator $T:E\to F$ is said to be linear if $T(\alpha x + \beta y)=\alpha T(x) + \beta T(y)$, for all
$x, y\in E$ and $\alpha, \beta \in \mathbb{R}$. Also, it is called positive if $x\ge 0$ implies $T(x)\ge 0$.

(ii) $T$ is called sublinear if $T(x+y)\le T(x)+T(y)$ (subadditive) and $T(\lambda x)=\lambda T(x)$ (positively homogeneous), for all
$x, y\in E$ and $\lambda \ge 0$. Also, it is called monotone if $x\ge y$ implies $T(x)\ge T(y)$.

Evidently that for the definition of monotonicity and sublinearity of $T$, is good enough to suppose that $E$ and $F$ are positive cones in
ordered vector spaces.

{\bf Remark 2.1.} Any positive linear operator is monotone sublinear but the converse is not true. Many concrete examples of monotone sublinear operators which are not positive linear are presented in, e.g.,  \cite{Gal-Nic-Med}--\cite{submit}. Also,
if $T$ is positively homogenous only, then $T(0)=0$ and for linear operators, the positivity is equivalent to monotonicity.

\section{Extension properties}

Firstly, for our purpose, we recall that a Riesz space (or a vector lattice) is an ordered vector space $E$ with
the additional property that for each pair of vectors $x, y \in E$, the supremum
and the infimum of the set $\{x, y\}$ both exist in $E$. In a Riesz space, two elements $x$ and $y$ are said to be disjoint (in
symbols $x \bot y$) whenever $|x| \wedge |y| = 0$ holds. A Riesz space is called Dedekind complete whenever every nonempty
bounded above subset has a supremum (or, equivalently, whenever every
nonempty bounded below subset has an infimum). Also, a Riesz space is said to be Dedekind $\sigma$-complete if every countable
subset that is bounded above has a supremum.

For details concerning linear ordered spaces see, e.g.,  Chapter 1, Section 1.1 in the book of Aliprantis-Burkinshaw \cite{Book-Al}.

We start with the following generalization of the fundamental
extension theorem of L. V. Kantorovich in \cite{Kant-1} (see also Aliprantis-Burkinshaw \cite{Book-Al}, Theorem 1.10, p. 9) proved there for positive linear operators.

{\bf Theorem 3.1.}
{\it Suppose that $E$ and $F$ are two Riesz spaces. Assume also that $T : E^{+} \to F^{+}$ is a monotone sublinear mapping.
Then, $S(x) = T(x^{+}), x \in  E$ is a monotone sublinear extension of $T$ on $E$.}

{\bf Proof.} Since $(x+y)^{+}\le x^{+}+y^{+}$, by the monotonicity of $T$ it follows
$T((x+y)^{+})\le T(x^{+})+T(y^{+})$, which immediately leads to
$$S(x+y)\le S(x)+S(y),$$
that is $S$ is subadditive on $E$. Also, for $\alpha \ge 0$, from the positive homogeneity of $T$ it easily follows that
$S$ is positive homogeneous too.

We prove now that $S$ is monotone on $E$. Thus, let $x, y\in E$ with $y\le x$. It follows $y^{+}=y\vee 0\le x\vee 0=x^{+}$,
which implies $S(y)=T(y^{+})\le T(x^{+})=S(x)$.

Evidently that $S(x)=T(x)$ for all $x\in E^{+}$ and $S(x)\ge 0$, for all $x\in E$.
$\hfill \square$

{\bf Remark 3.1.} Notice that in the case of positive linear operators in Aliprantis-Burkinshaw \cite{Book-Al}, Theorem 1.10, p. 9, the proof is completely different, where the extension is given by $S(x)=S(x^{+})-S(x^{-})$.

As an application of Theorem 3.1, the next result presents an interesting local approximation property of
monotone and sublinear operators.

{\bf Theorem 3.2.}  {\it Let $T : E \to F$ be a monotone sublinear operator between two Riesz
spaces with $F$ Dedekind $\sigma$-complete. Then for each $x\in E^{+}$ there exists a
monotone sublinear operator $S : E \to F$ depending on $x$, such that
$$0\le S(y) \le  T(y^{+}) \mbox{ for all } y\in E; \quad \quad  S(x) = T(x); \quad \quad \quad S(y) = 0 \mbox{ for all } y  \bot x.$$}

{\bf Proof.} Let $x \in E^{+}$ be fixed. As in the proof of Theorem 1.22, pp. 19-20 in Aliprantis-Burkinshaw \cite{Book-Al}, let us define
$S : E^{+} \to F^{+}$ defined for fixed  $x \in E^{+}$  by $S(y) = \sup\{T(y \wedge n x): n\in \mathbb{N} \}$.

(The supremum exists since $F$ is Dedekind $\sigma$-complete and the sequence
${T(y \wedge n x)}$ is bounded above in $F$ by $T(y)$.)

We claim that $S$ is subadditive.
To see this, let $y, z \in E^{+}$. From $(y + z) \wedge n x \le y \wedge n x + z \wedge n x$, we get
$$T((y + z) \wedge n x)\le T(y \wedge n x) + T(z \wedge n x),$$
which by passing to supremum after $n\in \mathbb{N}$ immediately implies $S(y + z) \le S(y) + S(z)$.

Then, for $\lambda \ge 0$ we get
$$S(\lambda y)=\sup\{T(\lambda y \wedge n x); n\in \mathbb{N}\}=\lambda\sup\{T(y+\frac{n}{\lambda} x); n\in \mathbb{N} \}.$$
But for each $n\in \mathbb{N}$, there exits $m\in \mathbb{N}$ with $\frac{n}{\lambda}\le m$, which implies $y+\frac{n}{\lambda} x\le y+m x$, $T(y+\frac{n}{\lambda} x)\le T(y+m x)$, $\sup\{T(y+\frac{n}{\lambda} x); n\in \mathbb{N} \}=\sup\{T(y+ m x); m\in \mathbb{N} \}$ and therefore
we immediately obtain $S(\lambda y)=\lambda S(y)$.

Now, if $y_1\le y_2$, by the monotonicity of $T$ it easily follows $S(y_1)\le S(y_2)$, that is $S$ is monotone too.

By Theorem 3.1 the mapping $S$ extends on $E$ to a monotone
sublinear operator.

Finally, it remains to prove the properties of $S$ from the end of the statement. Indeed, from the formula of definition for $S$, since
$y\wedge nx\le y$ it easily follows that $0\le S\le T$ for all $y\in E^{+}$. Also, by $x\wedge n x = x$ for all $n\in \mathbb{N}$, it immediately follows that $S(x)=T(x)$. Also, if $y\bot x$, it follows $y\wedge n x=0$ and therefore $S(y)=0$.
$\hfill \square$

{\bf Remark 3.2.} Theorem 3.2 is a generalization of Theorem 1.22 in \cite{Book-Al} valid for positive linear operators.

In what follows we prove for sublinear operators a variant of the classical Hahn-Banach theorem for linear operators.

{\bf Theorem 3.3.}
{\it Let $X$ be a real vector space, $Y\subset X$ a vector subspace of $X$, $F$ a Dedekind
complete Riesz space, $S:Y\to F$ sublinear operator on $Y$ and  $p: X\to  F$ a sublinear operator on $X$, such that $S(u)\le p(u)$, for all $u\in Y$.
Then there exists $T:X\to F$, such that $T(y)=S(y)$ for all $y\in Y$, $T$ is sublinear on $X$ and $T(x)\le p(x)$, for all $x\in X$.}

{\bf Proof.} The proof is a bit different from that in the linear case. Thus, for a fixed $x\in X\setminus Y$, define $Z$ as direct sum of $Y$ and $\mathbb{R}_{-} x$, that is $Z=\{z=y+\lambda x; y\in Y, \lambda \le 0\}$. Clearly that $Z$ is a cone in $X$.

Firstly we prove that there exists $T:Z\to F$, such that $T$ is sublinear on $Z$ and
$T(z)\le p(z)$, for all $z\in Z$. For this purpose, for any $u \in Y$, by the hypothesis we get
$S(u)\le p(u)\le p(u-x)+p(x)$,  which implies $S(u)-p(u-x)\le p(x), u\in Y$.

Therefore, taking here a supremum of the left-hand side over
$u$ denoted by $m$, we find that $m \le p(x)$. For a fixed $c\in [m, p(x)]$ define $T(z)=T(y+\lambda x)=S(y)+\lambda c$, $\lambda \le 0$.

Evidently that $Y\subset Z$ and $T(y)=S(y)$ for all $y\in Y$. Then,
$$T(z_1+z_2)=T(y_1+y_2+(\lambda_1+\lambda_2)x)=S(y_1 + y_2)+(\lambda_1+\lambda_2)c$$
$$\le S(y_1)+S(y_2)+\lambda_1 c + \lambda_2 c = T(z_1)+T(z_2),$$
which shows that $T$ is subadditive on $Z$.

Also, for any $\alpha\in \mathbb{R}_{+}$ and $z=y+\lambda x$ with $\lambda \le 0$, we get
$$T(\alpha z)=T(\alpha y+\alpha \lambda x)=S(\alpha y)+\alpha \lambda c=\alpha (S(y)+\lambda c)=\alpha T(z),$$
which shows that $T$ is positively homogeneous on $Z$.

In what follows we show that $T(z)\le p(z)$ for all $z=u+\lambda x\in Z$ with  $\lambda < 0$.

By the definition of $c$, we get $S(u)-p(u-x)\le c \le p(x), \mbox{ for all } u\in Y$.
For $\lambda < 0$, replacing in the left-hand side  $u$ by $-\lambda^{-1} u$, we get
$$c\ge S(-\lambda^{-1} u)-p(-\lambda^{-1}(u+\lambda x))=-\lambda^{-1}(S(u)-p(u+\lambda x)),$$
which immediately implies that $T(z)\le p(z)$.

Finally, a standard application of Zorn's lemma as in the classical linear case, guarantees the existence
of an extension of $S$ to $X$ with the desired properties. Indeed, denote by $\mathcal{F}$ the set of all $T:Z\to F$; $Y$ subspace of $Z$, $T$ is sublinear on  $Z$, $T(z)=S(z), z\in Y$, $T(z)\le p(z), z\in Z$.

If we consider on $\mathcal{F}$ the order
$T_1, T_2 \in \mathcal{F}$, $T_1 \le T_2$ if and only if $T_2$ is an extension of $T_1$, then $(\mathcal{F}, \le)$ is partially ordered and by Zorn's lemma it has a maximal element, let us denote it by $T^{*}:Y^{*}\to F$. We have here $Y^{*}=X$, because if we suppose that there exists $x\in X\setminus Y^{*}$, by the previous reasonings there exists $T^{**}:Z^{*}\to F$, with $Z^{*}=Y^{*}+\mathbb{R}_{-} x$ such that $T^{**}(z)=T^{*}(z)$, for all $z\in Z^{*}$. But since $T^{**}\in \mathcal{F}$, it follows that $T^{*}\le T^{**}$, which by the maximality of $T^{*}$ implies in fact that
$T^{*}=T^{**}$ and therefore $Y^{*}=Z^{*}$, i.e. $x\in Y^{*}$, which is a contradiction.
$\hfill \square$

{\bf Remark 3.3.}
If $S$ is  linear, this results was obtained by Theorem 1.25 in \cite{Book-Al}.

As a consequence of Theorem 3.3, we get the following Hahn-Banach type theorem for monotone/positive sublinear operators.

{\bf Theorem 3.4.}
{\it Let $T : E \to F$ be a positive sublinear operator between two Riesz
spaces with $F$ Dedekind complete. Assume also that $G$ is a Riesz subspace
of $E$ and that $S : G \to F$ is a monotone sublinear operator satisfying $0 \le  S(x) \le T(x)$, for all
$x$ in $G^{+}$. Then $S$ can be extended to a positive sublinear operator $Q:E\to F$ such
that $0 \le Q(x) \le T(x)$ holds for all $x\in E^{+}$.}

{\bf Proof.} As in the linear case, define $p: E \to F$ by $p(x) = T(x^{+})$ and note firstly that $p$ is sublinear. Indeed, this is immediate from $(x+y)^{+}\le x^{+}+
y^{+}$ and $(\lambda x)^{+}=\lambda x^{+}$, for all $x, y\in E$ and $\lambda \ge 0$. Also, since $x\le x^{+}$, for all $x\in G$ we have
$$S(x)=S(x^{+}-x^{-})\le S(x^{+})\le T(x^{+})=p(x).$$
Therefore, by Theorem 3.3, there exists a
sublinear extension of $S$ to all of $E$ (which we denote by $Q$ ) satisfying
$Q(x) \le  p(x)$,  for all $x \in  E$.

By $0=Q(x-x)\le Q(x)+Q(-x)$, it follows that $-Q(x)\le Q(-x)$  and for $x\in E^{+}$ we get
$$-Q(x) \le Q(-x) \le p(-x) = T((-x)^{+})= T(0) = 0$$
and so $0 \le Q(x) \le p(x) = T(x)$ holds for all $x\in E^{+}$, as desired.
$\hfill \square$

{\bf Remark 3.4.} Theorem 3.4 is a generalization of Theorem 1.26 in \cite{Book-Al}, p. 24.

Recall now that a subset $A$ of a Riesz space $E$ is said to be order bounded if there exists $m, M\in E$ such that
$m\le x \le M$ for all $x\in A$. Also, an operator $T : E \to F$ between two Riesz spaces is said
to be order bounded if it maps order bounded subsets of $E$ to order bounded
subsets of $F$.

{\bf Theorem 3.5.} {\it Let $E$ and $F$ be Riesz spaces with $F$ Dedekind complete,
$G$ a Riesz subspace of $E$ and $T : G \to F$ be a monotone sublinear operator.  Then for the
following statements

(1) $T$ extends to a monotone sublinear operator from  $E$ to $F$.

(2) $T$ extends to an order bounded sublinear monotone operator $S$ from  $E$ to  $F$.

(3) There exists a monotone sublinear mapping  $p: E \to F$ satisfying
$T(x) \le p(x)$ for all  $x \in  G.$

(4) $T$ extends to a positive sublinear operator from $E$ to $F$,

the following implications hold :

$(1) \Longrightarrow  (2) \Longrightarrow (3) \Longrightarrow (4)$.}

{\bf Proof.} $(1) \Longrightarrow  (2)$. Is obvious.

$(2) \Longrightarrow  (3)$. Let $S$ be an order bounded sublinear monotone operator satisfying $S(x) = T(x)$ for all $x \in G$.  Then
the mapping $p: E \to F$ defined by $p(x) = |S(x^{+})|$ is monotone, sublinear
and satisfies
$$T(x) \le  T(x^{+}) = S(x^{+}) \le |S(x^{+})| = p(x),$$
for all $x \in  G$.

$(3) \Longrightarrow (4)$.  Let $p: E \to F$ be a monotone sublinear mapping satisfying
$T(x) \le  p(x)$ for all $x \in G$. Then the formula $q(x) = p(x^{+})$ defines a
sublinear mapping from $E$ to $F$ such that we have
$$T(x) \le T(x^{+}) \le p(x^{+}) = q(x),$$
for all $x \in  G$. Thus, by the Hahn-Banach extension Theorem 3.3,
there exists a sublinear extension $R$ of $T$ on $E$,  satisfying $R(x) \le q(x)$, for all
$x \in E$. In particular, if $x ∈\in E^{+}$, then by the relation $0=R(x-x)\le R(x)+R(-x)$, we get
$$-R(x) \le R(-x) \le q(-x) = (p(-x)^{+})=p(0)=0,$$
which implies $R(x) \ge 0$. That is, $R$ is a positive sublinear extension of $T$ on $E$,
and the proof is finished. $\hfill \square$

{\bf Remark 3.5.}
If in the above theorem $T$ is a positive linear operators, then we recapture
Theorem 1.27 in Aliprantis-Burkinshaw \cite{Book-Al}.

A subset $A$ of a Riesz space is called solid whenever $|x| \le |y|$ and $y \in A$
imply $x \in A$. A solid vector subspace of a Riesz space is called
an ideal. From the lattice identity $x \vee y = \frac{1}{2} (x + y + |x - y|)$, it follows
immediately that every ideal is a Riesz subspace.

The next result deals with restrictions of monotone sublinear operators to ideals.

{\bf Theorem 3.6.} {\it If $T : E \to F$ is a monotone sublinear operator between two Riesz
spaces with $F$ Dedekind complete, then for every ideal $A$ of $E$ the formula
$$T_{A}(x) = \sup \{T(y): y \in A \mbox{ and } 0 \le y \le x\}, x\in E^{+},$$
defines (by its extension on $E$ denoted also by $T_{A}$ and proved to exist by Theorem 3.1)  a monotone sublinear operator from $E$ to $F$. Moreover, we have :

(a) $0 \le T_{A}(y) \le T(y^{+})$, for all $y\in E$.

(b) $T_{A} = T$ on $A$ and $T_{A} = 0$ on $A^{d} := \{x \in E; x \bot y \mbox{ for all } y \in A\}$.

(c) If $B$ is another ideal with $A \subset  B$, then $T_{A} \le T_{B}$ holds.}

{\bf Proof.} According to Theorem 3.1,  we need to show that $T_{A}$
is monotone and sublinear on $E^{+}$. The monotonicity immediately follows from the definition of $T_{A}$.

It remains to prove the sublinearity on $E^{+}$.  For this purpose, firstly it is easy to see that
$$T_{A}(x) = \sup\{T(x  \wedge y); y ∈\in A^{+}\}$$
holds for all $x \in E^{+}$.

Let $x, y \in E^{+}$. If $z \in A^{+}$, then the inequality
$$(x + y) \wedge  z \le x \wedge z +  y \wedge z$$
implies that
$$T((x + y) \wedge z) \le T(x \wedge z) + T(y  \wedge z) \le T_{A}(x) + T_{A}(y),$$
and hence
$$T_{A}(x + y) \le T_{A}(x) + T_{A}(y),$$
that is the subadditivity.

Now, since for $\lambda >0$ and $y\in A^{+}$ we have $\frac{y}{\lambda}\in A^{+}$, we immediately get
$$T_{A}(\lambda x)=\sup\{T(\lambda x \wedge y) ; y\in A^{+}\}=\lambda \sup\{T(x \wedge \frac{y}{\lambda}); y\in A^{+}\}$$
$$=\lambda \sup\{T(x \wedge y); y\in A^{+}\}=\lambda T_{A}(x),$$
which shows that $T_{A}$ is positive homogeneous too.

Now, the properties $(a)-(b)$ are easy consequences of Theorem 3.2, while the property (c) is an easy consequence of the formula defining
$T_{A}$. $\hfill \square$

{\bf Remark 3.6.} Theorem 3.6 generalizes Theorem 1.28, pp. 25-26 in \cite{Book-Al}.

In what follows, one says that a vector subspace $G$ of an ordered vector space $E$ is
majorizing $E$, if for each $x \in E$, there exists some $y \in  G$ with $x \le y$
(or, equivalently, if for each $x \in  E$, there exists some $y \in G$ with $y \le x$).

{\bf Theorem 3.7.} {\it Let $E$ and $F$ be two ordered vector spaces
with $F$ a Dedekind complete Riesz space. If $G$ is a majorizing vector subspace
of $E$ and $T : G \to F$ is a monotone sublinear operator, then $T$ has a positive sublinear
extension to all of $E$.}

{\bf Proof.}  Fix $x \in  E$ and let $y \in G$ satisfy $x \le y$. Since $G$ is majorizing, there
exists a vector $u \in G$ with $u \le x$. Hence, $u \le y$ and the monotonicity of $T$
implies $T(u) \le T(y)$ for all $y \in G$ with $x \le y$. In particular, it follows that
$\inf\{T(y): y \in G \mbox{ and } x \le y\}$
exists in $F$ for each $x \in  E$ and we can define the mapping
$p: E \to F$ by the formula
$$p(x) = \inf\{T(y): y \in G \mbox{ and } x \le y\}, x\in E.$$

For each $x \in  G$ we have  $p(x) = T(x)$. Then, for $\lambda > 0$, we get
$$p(\lambda x)=\inf\{T(y); y \in G \mbox{ and } \lambda x\le y\}=\inf\{T(\lambda \cdot y/\lambda); y \in G \mbox{ and } x\le y/\lambda\}$$
$$=\lambda \inf\{T(y/\lambda); y \in G \mbox{ and } x\le y/\lambda\}=\lambda p(x).$$
Also, let $x_{1}$ and $x_{2}$ be arbitrary in $E$. Then for any $y_{1}, y_{2}\in G$
with $x_{1}\leq y_{1}$ and $x_{2}\leq y_{2}$, we have $x_{1}+x_{2}\leq
y_{1}+y_{2}$ and therefore
\[
p\left( x_{1}+x_{2}\right) \leq T(y_{1}+y_{2})\leq T(y_{1})+T(y_{2}).
\]%
Passing  to infimum successively with $y_{1}$ and $y_{2}$, it follows that
$p\left( x_{1}+x_{2}\right) \leq p(x_{1})+p(x_{2}),$ i.e., $p$ is sublinear.

Now, by Theorem 3.3, the operator $T$ has
a sublinear extension $S$ on $E$, satisfying $S(z) \le  p(z)$ for all $z \in  E$. If
$z \in E^{+}$, then $-z \le 0$, and so from $0=S(z-z)\le S(z)+S(-z)$ it follows
$$-S(z) \le S(-z) \le p(-z) \le T(0) = 0,$$
and we see that $S(z)\ge  0$. This shows that $S$ is a positive sublinear extension of $T$
on $E$.  $\hfill \square$

{\bf Remark 3.7.} Theorem 3.7 generalizes Theorem 1.32 in \cite{Book-Al} (see also Kantorovich \cite{Kant-2}).

Now, for a monotone sublinear operator $T : G \to F$, where $G$ is a vector subspace
of an ordered vector space $E$ and $F$ is a Dedekind complete Riesz
space, let us denote by  $\mathcal{E}_{sub}(T)$, the collection of all monotone and sublinear extensions of $T$
on $E$. Namely,
$$\mathcal{E}_{sub}(T) := \{S \mbox{ monotone, sublinear and } S = T \mbox{ on } G\}.$$

Clearly the set $\mathcal{E}_{sub}(T)$ is a convex subset, that is $\lambda S + (1 - \lambda)R \in \mathcal{E}_{sub}(T)$
for all $S, R \in \mathcal{E}_{sub}(T)$ and all $\lambda \in [0, 1]$.

In what follows we prove that any extendable monotone sublinear operator whose domain is an
ideal has a smallest extension.

{\bf Theorem 3.8.} {\it Let $E$ and $F$ be two Riesz spaces with $F$ Dedekind complete,
let $A$ be an ideal of $E$, and let $T : A \to F$ be a monotone sublinear operator. If $\mathcal{E}_{sub}(T) \not= \emptyset$,
then $T$ has a smallest extension. Moreover, if in this case $S = \min \mathcal{E}_{sub}(T)$,
then
$$S(x) = \sup\{T(y) : y \in A \mbox{ and } 0 \le y \le x\}, \mbox{ for all } x \in E^{+}.$$}

{\bf Proof.} Since $\mathcal{E}_{sub}(T) \not= \emptyset$, $T$ has (at least) one monotone and sublinear  extension, the formula
$$T_{A}(x) = \sup\{T(y) : y \in A \mbox{ and } 0 \le y \le x\},$$
$x \in  E^{+}$ ,
by Theorem 3.6 extends to a monotone sublinear operator from $E$ to $F$ (denoted here also by $T_{A}$), satisfying $T_{A} = T$ on $A$, and so
$T_{A}\in  \mathcal{E}_{sub}(T)$.
Now if $S \in  \mathcal{E}_{sub}(T)$, then $S = T$ holds on $A$, and hence $T_{A} = S_{A} \le S$.
Therefore, $T_{A} = \min \mathcal{E}_{sub}(T)$ holds, as desired. $\hfill \square$

{\bf Remark 3.8.} Theorem 3.8 generalizes Theorem 1.30, pp. 27 in \cite{Book-Al}.

In what follows, recall that a vector $e$ of a convex set $C$ is called an extreme point of
$C$ if $e = \lambda x +(1 - \lambda)y$ with $x, y \in C$ and $0 < \lambda < 1$
implies $x = y = e$.

{\bf Theorem 3.9.} {\it Let $E$ and $F$ be two Riesz
spaces with $F$  Dedekind complete, $G$ a vector subspace of $E$ and
$T : G \to F$ be a monotone sublinear operator.

If $S\in \mathcal{E}_{sub}(T)$, satisfies the condition that for all $x \in E$ we have $\inf\{S(|x - y|); y \in G\} = 0$, then $S$ is an extreme point of $\mathcal{E}_{sub}(T)$.}

{\bf Proof.} Let $S\in \mathcal{E}_{sub}(T)$ satisfying the hypothesis and assume that $S = \lambda Q + (1 - \lambda) R$ with
$Q, R \in \mathcal{E}_{sub}(T)$ and $0 < \lambda < 1$. Then for each $x, y \in  E$ we have
$$|Q(x) - Q(y)|\le Q(|x - y|) = (\frac{1}{\lambda}S - \frac{1-\lambda}{\lambda}R)(|x - y|) \le \frac{1}{\lambda}S(|x - y|).$$
Here the inequality $|Q(x) - Q(y)|\le Q(|x - y|)$ follows immediately applying the monotonicity of $Q$ to both inequalities $x-y\le |x-y|$ and $y-x\le |x-y|$.

In particular, if $x \in E$ and $y \in G$, then from $S(y) = Q(y) = T(y)$ it follows
that
$$|S(x) - Q(x)|\le |S(x) - S(y)|+|Q(y) - Q(x)|\le (1 + 1/\lambda)S(|x - y|).$$
Taking into account our hypothesis, the last inequality yields $S(x) = Q(x)$
for each $x \in E$, and this shows that $S$ is an extreme point of $\mathcal{E}_{sub}(T)$. $\hfill \square$

{\bf Remark 3.9.} For positive linear operators, the converse of the statement holds too, see Lipecki-Plachky-Thomsen \cite{Lip} (see also Theorem 1.31, p. 27 in \cite{Book-Al}).

{\bf Acknowledgement.} The author thanks Constantin P. Niculescu for an useful discussion on the proof of Theorem 3.7.




\begin{thebibliography}{99}                                                                                               %

\bibitem{Book-Al}Aliprantis, C.D., Burkinshaw, O.: Positive Operators. Springer, Dordrecht (2006).

\bibitem {Gal-Nic-Med}Gal, S.G., Niculescu, C.P.: A nonlinear extension of
Korovkin's theorem. Mediterr. J. Math. \textbf{17}, Article no. 145 (2020), 14 pages.

\bibitem {Gal-Nic-JMAA}Gal, S.G., Niculescu, C.P.: Choquet operators
associated to vector capacities. J. Math. Anal. Appl. \textbf{500}, article
no. 125153 (2021), 24 pages.

\bibitem {Gal-Nic-Aeq}Gal, S.G., Niculescu, C.P.: A note on the Choquet type
operators. Aequationes Math. \textbf{95}, 433--447 (2021).

\bibitem {Gal-Nic-RACSAM}Gal, S.G., Niculescu, C.P.: Nonlinear versions of
Korovkin's abstract theorems. Rev. Real Acad. Cienc. Exactas Fis. Nat. Ser. A-Mat.
\textbf{116}, Article number 68 (2022)

\bibitem{submit}Gal, S.G., Niculescu, C.P.: Korovkin type theorems for weakly nonlinear and monotone operators.
Mediterr. J. Math. {\bf 20}, Article no. 56 (2023), 20 pages .

\bibitem{Kant-1}Kantorovich, L.V.: Linear operators in semi-ordered spaces, Math. Sbornik,
209--284 (1940)

\bibitem{Kant-2}Kantorovich, L.V.: On the moment problem for a finite interval, Dokl. Akad. Nauk
SSSR, {\bf 14}, 531--537 (1937)

\bibitem{Lip} Lipecki, Z., Plachky, D., Thomsen, W: Extension of positive operators and extreme
points I. Colloq. Math. {\bf 42}, 279--284 (1979).






\end{thebibliography}
\end{document}